\documentclass[12pt, pra, showpacs, longbibliography]{revtex4-1}
\usepackage{epsfig}
\usepackage{amsmath}
\usepackage{dcolumn}
\usepackage{bm}
\usepackage{graphicx}
\usepackage{wrapfig}
\usepackage{graphicx}
\usepackage{epstopdf}
\usepackage{epsfig}
\usepackage{calrsfs}
\usepackage{IEEEtrantools}
\usepackage{amssymb}
\usepackage{url} 
\usepackage{hyperref}

\usepackage{epsfig}
\usepackage{amsmath}
\usepackage{dcolumn}
\usepackage{bm}
\usepackage{graphicx}
\usepackage{wrapfig}
\usepackage{graphicx}
\usepackage{epstopdf}
\usepackage{epsfig}
\usepackage{calrsfs}
\usepackage{IEEEtrantools}
\usepackage{graphicx}
\usepackage{caption}
\usepackage{mathtools}
\usepackage{dsfont}
\usepackage{subfigure}
\graphicspath{{"/Users/jameswatson/Documents/My Documents/Gavrila-LC/MathPlots/"}}

\newcommand{\nn}{\nonumber}

\newcommand{\Sch}{Schr\"{o}dinger}
\newcommand{\GL}{Gr\"{u}nwald–-Letnikov}
\newcommand{\RL}{Riemann--Liouville}
\newcommand{\Cp}{Caputo}

\newcommand{\Cy}{Cauchy}
\newcommand{\Fr}{Fourier}

\newcommand{\cl}{\colon}

\begin{document}



\title{Fractional derivative of composite functions: \\ exact results and physical applications}

\author{Gavriil Shchedrin, Nathanael C. Smith, Anastasia Gladkina, and Lincoln D. Carr}
\affiliation{Colorado School of Mines, Golden, Colorado 80401, USA}

\begin{abstract}
We examine the fractional derivative of composite functions and present 
a generalization of the product and chain rules for the {\Cp} fractional derivative. 
These results are particularly important for a comprehensive description of transport phenomena through multiscale physical systems and biological structures, e.g., porous materials, disordered media, and clusters of neurons. In order to address the multiscale dynamics that occur in these systems, we derive product and chain rules for the {\Cp} fractional derivative of composite functions. Unlike the Leibniz and di Bruno formulae that characterize an integer-order derivative of a product of functions and composite functions, respectively, and which result in a finite series of lower order derivatives, the fractional analogs of these formulae produce an infinite series of fractional derivatives of the constituent functions. 
We demonstrate the obtained results by the exact evaluation of the {\Cp} fractional derivative of hyperbolic tangent, which describes dark soliton propagation in non-linear media.
\end{abstract}

\maketitle


\clearpage

\section{Introduction}

Transport through multi-scale physical and biological systems,
such as tissues \cite{murray1926physiological, west1997general, sherman1981connecting}, clusters of neurons \cite{lundstrom2008fractional, anastasio1994fractional}, porous materials \cite{aizenberg2005skeleton, liu2011bio, mcculloh2003water}, disordered media \cite{metzler2000random, saichev1997fractional},
and ultimately the Solar System \cite{mathai2009h} is governed by fractional partial differential equations (FPDEs) \cite{samko1993fractional, kilbas2006theory, mcbride1979fractional}. These vastly different physical systems share distinct emergent phenomena, including, but not limited to, long-range correlations, non-locality, fractional geometry, non-Gaussian statistics, and non-Fickian transport \cite{west2014colloquium, herrmann2014fractional}. The general framework of FPDEs provides a thorough account of these properties in a cohesive and self-consistent way \cite{samko1993fractional, kilbas2006theory, mcbride1979fractional, kiryakova1993generalized, herrmann2014fractional}. Moreover, this framework is not only capable of describing the properties of existing multi-scale physical systems, but also allows one to design and build advanced synthetic materials with {\it prescribed} physical properties, such as enhanced mass exchange and charge transfer rate \cite{liu2011bio, zheng2017bio}.

The development of the general framework of FPDEs 
has brought a rich variety of fractional derivatives -  from the discrete {\GL} fractional derivative defined in a coordinate space to a continuous {\Fr} fractional derivative defined in a frequency domain \cite{samko1993fractional, kilbas2006theory, mcbride1979fractional, herrmann2014fractional}. The non-local nature of the {\Cp} fractional derivative, combined with its convergence at the origin, makes it the most suitable choice of a fractional derivative in a wide range of physical applications in condensed matter, astrophysics, biophysics, and material science \cite{herrmann2014fractional, mathai2009h, west2014colloquium}. 
Thus, exact results for the {\Cp} fractional derivative play a key role in the description of  multi-scale physical and biological systems. A number of versatile and robust numerical techniques and analytical methods have been developed that allow one to evaluate a whole range of fractional derivatives for a wide spectrum of functions
\cite{kilbas2004generalized, srivastava2009fractional, shukla2007generalization, mainardi2000mittag, haubold2011mittag, saadatmandi2010new, shawagfeh2002analytical, momani2006analytical, odibat2006application, rossikhin2010application, daftardar2005adomian, dehghan2010solving, momani2005explicit, li2010stability, diethelm2005algorithms, agrawal2002solution, kiryakova2010special, kiryakova1997all, kiryakova1993generalized,gladkina2017expansion}. Recently, an analytic method has been developed for an exact evaluation of the {\Cp} fractional derivative 
based on the generalized Euler's integral transform (EIT) \cite{shchedrin2017exact}. Despite the fact that this method led to the exact evaluation of the {\Cp} fractional derivative of a broad class of elementary functions, such as Gaussian, quartic Gaussian, Lorentzian, and hyperbolic functions, it was not a truly universal method. Indeed, this method was limited to a class of functions that can be expressed in terms of a hypergeometric function with a {\it power-law argument}. 
While, individually, elementary functions can be represented in terms of a single hypergeometric function with a power-law argument, their combination, in the most general case, cannot be brought to such a form. These combinations are especially important 
for a number of physical applications, e.g. Gaussian wavepacket propagation described by the fractional {\Sch} equation \cite{herrmann2014fractional, laskin2002fractional}. 
If we assume that the amplitude is a slowly varying function compared to the phase of the wavepacket, one can apply the slowly varying envelope approximation (SVEA) \cite{landau1965quantum}. In this approximation, one can decouple the amplitude from the highly oscillatory phase that allows one to solve the fractional {\Sch} equation \cite{herrmann2014fractional}. The SVEA method, along with the quasi-classical approximation and the method of separation of variables, are just a few examples among the myriad of methods for solving FPDEs that rely on the decomposition of the wavefunction in terms of a product of the trial functions. Since the EIT method cannot be used in the evaluation of the {\Cp} fractional derivative of composite functions, we address this problem by deriving both product and chain rules for the {\Cp}  fractional derivative. We shall point out that a form of the product rule has been obtained previously \cite{samko1993fractional, kilbas2006theory}. However, practical implementation 
of the fractional chain and product rules was limited due to the fact that the results
were expressed in terms of a sum of the repeated integral of the generalized hypergeometric function with a power-law argument. In this paper, we solve this problem by applying the EIT method \cite{shchedrin2017exact} to the repeated integral that results in a single hypergeometric function of a higher order. As a result, we are able to implement the fractional product and chain rules in the application of the {\Cp} fractional derivative in a much simpler and more practical way. Based on the obtained results, we are able to extend the applicability of the generalized EIT to a domain of functions that cannot be expressed in terms of a single 
hypergeometric function with a power-law argument.

The rest of this paper has the following structure. In sections {\bf \ref{section_product}} and {\bf \ref{section_chain}} we expand the Caputo fractional derivative into an infinite series of integer order derivatives and derive the fractional product and chain rules, correspondingly. In section {\bf \ref{section_appl}} we apply the EIT method to the exact evaluation of the repeated integral of the generalized hypergeometric function with a power-law argument, which enables the practical implementation of both fractional product and chain rules. Finally, we demonstrate the fractional product and chain rules by evaluating the {\Cp} fractional derivative of hyperbolic tangent. In section {\bf \ref{section_concl}} we summarize the obtained results.

\section{The fractional product rule}\label{section_product}

In this paper we will focus on the {\Cp} fractional derivative due to its fundamental role in physical applications \cite{herrmann2014fractional}. The {\Cp} fractional derivative of fractional order $0<\alpha<{1}$ for $x>0$ is defined as \cite{samko1993fractional, kilbas2006theory, herrmann2014fractional},
\begin{eqnarray}\label{capdefinition}
{}^{\rm C}D^{\alpha}_{x}f(x) = \frac{1}{\Gamma(1-\alpha)}\int^{x}_{0}{dt}\;(x-t)^{-\alpha}\frac{df(t)}{dt} 
.
\end{eqnarray}
To extend this definition to the entire domain of $x\in (-\infty,\infty)$, a factor of $\mathrm{Re}[\exp(-i\alpha)]$ is applied to the result of the derivative in the $x<0$ region, but for simplicity we will restrict ourselves to the definition in Eq.\eqref{capdefinition} \citep{samko1993fractional}. We shall point out that the {\Cp} fractional derivative is defined only for the non-integer values of fractional order, i.e., $\alpha\notin{\mathds{N}}$ \cite{kilbas2006theory, samko1993fractional, herrmann2014fractional}. For integer values of order $\alpha=n\in{\mathds{N}}$, the {\Cp} fractional derivative is defined 
such that it matches exactly with the integer derivative of the $n^{\rm th}$ order \cite{samko1993fractional}. Previously we found that the {\Cp} fractional derivative 
of the order $\alpha=n\in{\mathds{N}}$ is consistent with the integer $n^{\rm th}$ order derivative, except for the case $\alpha=0$ \cite{shchedrin2017exact}. In the special case of $\alpha=0$, the {\Cp} fractional derivative differs from the function $f(x)$ by its value at the origin, i.e., ${}^{\rm C}D^{\alpha=0}_{x}f(x) = f(x)-f(0)$.
The extra requirement for the {\Cp} fractional derivative of integer order ${}^{\rm C}D^{\alpha=n}_{x}f(x) \equiv{}d^{n}f(x)/dx^{n}$ leads to a discontinuity in the vicinity of the parameter $\alpha=0$. However,  the {\Cp} fractional derivative defined by 
Eq.~(\ref{capdefinition}) adopted for both integer and non-integer values of $0\leq{}\alpha\leq{1}$ leads to a continuous transformation of the fractional derivative between ${}^{\rm C}D^{\alpha=0}_{x}f(x)$ and ${}^{\rm C}D^{\alpha=1}_{x}f(x)$.  In order to avoid the discontinuity in the vicinity of $\alpha=0$, we apply the definition of the {\Cp} fractional derivative given by Eq.~(\ref{capdefinition}) for all values of $0\leq{}\alpha\leq{}{1}$. 

Previously we established that a wide range of fractional derivatives - from discrete {\GL} to the continuous {\RL} and {\Cp} fractional derivatives \cite{kilbas2006theory, samko1993fractional, herrmann2014fractional} -  
can be equivalently expressed in terms of an infinite expansion of integer order derivatives \cite{gladkina2017expansion}. The expansion for the Caputo derivative of order $0\leq{}\alpha\leq{}{1}$ takes the form 
%
\begin{eqnarray}\label{frac_integer_sum1}
D_x^{\alpha}[f(x)] = 
\sum^{\infty}_{k=0} 
\frac{ \sin [\pi  (\alpha-k)]}{\pi  (\alpha-k)} \left(\frac{ \Gamma (\alpha+1)}{\Gamma (k+1)} 
\right)
x^{k-\alpha} 
\frac{d^{k}}{d x^{k}}
[f(x) - f(0)],
\end{eqnarray}
where, for convenience and to be consistent with the adoption of this derivative for $\alpha=0$, we will redefine $f(x)\equiv f(x)-f(0)$ for the duration of this paper. We shall point out that, by using this expansion, one can derive the fractional product rule for a product of an arbitrary number of functions. However, for the sake of simplicity 
and without loss of generality, we restrict ourselves to the Caputo fractional derivative of a product of two functions,
\begin{align}
D_x^{\alpha}[f(x)\cdot g(x)] &= 
\sum^{\infty}_{k= 0} 
\frac{ \sin [\pi  (\alpha-k)]}{\pi  (\alpha-k)}  \left(\frac{ \Gamma (\alpha +1)}{\Gamma (k+1)} 
\right)
x^{k-\alpha} 
\frac{d^{k}}{d x^{k}}
[f(x)\cdot g(x)] =\\\nn & = 
\sum^{\infty}_{k= 0} 
\frac{ \sin [\pi  (\alpha-k)]}{\pi  (\alpha-k)} \left(\frac{ \Gamma (\alpha+1)}{\Gamma (k+1)} 
\right)
x^{k-\alpha} 
\sum_{l=0}^{k} C^{l}_{k}f^{(k-l)}(x) g^{(l)}(x),
\end{align}
where $ C^{l}_{k}= \binom{k}{l} = {k!}/{[l! (k-l)!}]$ is the binomial coefficient. 
By exchanging the order in the summation,
\begin{eqnarray}
\sum^{\infty}_{k= 0} \sum_{l=0}^{k}  = 
\sum^{\infty}_{l= 0} \sum_{k=l}^{\infty}
, 
\end{eqnarray}
we obtain, 
\begin{eqnarray}
D_x^{\alpha}[f(x)\cdot g(x)] = 
\sum^{\infty}_{l= 0} 
g^{(l)}(x) 
\sum_{k=l}^{\infty} 
\frac{ \sin [\pi  (\alpha-k)]}{\pi  (\alpha-k)}  \left(\frac{ \Gamma (\alpha+1)}{\Gamma (k+1)} 
\right)
x^{k-\alpha} 
 C^{l}_{k}f^{(k-l)}(x) 
 .
\end{eqnarray}
Next, we perform a shift in the dummy summation index $k\to{k+l}$, which directly leads to 
\begin{eqnarray}
D_x^{\alpha}[f(x)\cdot g(x)] = 
\sum^{\infty}_{l= 0} 
g^{(l)}(x) 
\sum_{k=0}^{\infty} 
\frac{\sin [\pi  ((\alpha-l)-k)]}{\pi  ((\alpha-l)-k)}
 \frac{\Gamma(\alpha +1)}{\Gamma(k+ 1)\Gamma(l + 1)}
x^{k-(\alpha-l)} 
f^{(k)}(x)   
.
\end{eqnarray}
By rewriting the ratio of the Gamma functions,
\begin{eqnarray}
 \frac{\Gamma(\alpha +1)}{\Gamma(k+ 1)\Gamma(l + 1)} = 
  \frac{\Gamma(\alpha-l+1)}{\Gamma(k+1)}
  \frac{\Gamma(\alpha+1)}{\Gamma(\alpha-l + 1) \Gamma(l + 1)} = 
  \frac{\Gamma(\alpha-l+1)}{\Gamma(k+1)}
  C^{l}_{\alpha} 
  ,
\end{eqnarray}
we arrive at the fractional product rule,
\begin{align}\label{caputo_product}
D_x^{\alpha}[f(x)\cdot g(x)] &= 
\sum^{\infty}_{l= 0} 
  C^{l}_{\alpha}  
g^{(l)}(x) 
\sum_{k=0}^{\infty} 
\frac{\sin [\pi  ((\alpha-l)-k)]}{\pi  ((\alpha-l)-k)}
  \frac{\Gamma(\alpha-l+1)}{\Gamma(k+1)}
x^{k-(\alpha-l)} 
f^{(k)}(x)  =
\nn\\ &= 
\sum^{\infty}_{l= 0} 
  C^{l}_{\alpha}  
g^{(l)}(x) 
D_x^{(\alpha-l)}[f(x)]
,
\end{align}
where we have recognized that the inner sum in Eq.(\ref{caputo_product}) is nothing but the fractional derivative of the $(\alpha-l)^{\text th}$ order of the function $f(x)$ as it follows directly from Eq.(\ref{frac_integer_sum1}). Thus, the fractional derivative of the $\alpha^{\text th}$ order  of a product of two functions is given by an infinite series of a product of 
an $l^{\rm th}$ order integer derivative of the first function and an $(\alpha-l)^{\rm th}$ order fractional derivative of the second function \cite{samko1993fractional, kilbas2006theory}. 
We shall point out that one can bring the obtained fractional product rule given by Eq.~(\ref{caputo_product}) into a form in which the {\Cp} fractional derivative acts on both functions 
in a symmetric fashion, similarly to the Leibniz rule \cite{samko1993fractional}. However, in this case, the lower bound of summation in Eq.~(\ref{caputo_product}) 
of $l=0$ turns into $l\to{-\infty}$, i.e., the infinite sum goes over the entire range of indices \cite{samko1993fractional}. Since the {\Cp} fractional derivative of an elementary function is given in terms of the generalized hypergeometric function \cite{shchedrin2017exact}, the symmetric form of the fractional product rule will produce an infinite sum of a product of them.
Thus, the symmetric form of the fractional product rule, while being completely equivalent to the asymmetric expansion given by Eq.~(\ref{caputo_product}), results in a much more complicated expression for the {\Cp} fractional derivative. Therefore, for the sake of simplicity we will focus on the asymmetric form of the fractional product rule
as it appears in Eq.(\ref{caputo_product}). In the special case of an integer value of the parameter $\alpha=n\in{\mathds{N}}$, the infinite series in 
the fractional product rule given by Eq.(\ref{caputo_product}) becomes finite due to the
fact that the binomial coefficient vanishes, i.e. $C^{l}_{n}=0$, for integer values of index $l$ exceeding $n$, 
\begin{equation}\label{caputo_product2}
D_x^{\alpha=n}[f(x)\cdot g(x)] = 
\sum^{\infty}_{l= 0} 
  C^{l}_{n}  
g^{(l)}(x) 
D_x^{(n-l)}[f(x)] = 
\sum^{n}_{l= 0} 
  C^{l}_{n}  
g^{(l)}(x) 
f^{(n-l)}(x) = (f(x)\cdot g(x))^{(n)}.
\end{equation}
As a result, the fractional derivative of the integer order $\alpha=n\in{\mathds{N}}$ reduces to the well-known Leibniz rule \cite{samko1993fractional}.

\section{The fractional chain rule}\label{section_chain}

The goal of this section is to derive the fractional chain rule for the {\Cp} fractional derivative. We begin our derivation with di Bruno's formula for the $k^{\rm th}$ order integer derivative of a composite function $f(g(x))$ \cite{faa1855sullo, di1857note, jordan1965calculus, wolfram},
\begin{eqnarray}
\frac{d^{k}}{dx^{k}}
f(g(x))=
\sum _{m=0}^{k} \frac{1}{m!}\left(\sum _{j=0}^m (-1)^j \
{C^{j}_{m}}
g(x)^j
\frac{d^{k}}{dx^{k}}
g(x)^{m-j}\right) f^{(m)}(g(x))
.
\end{eqnarray}
Thus, the direct application of di Bruno's formula leads to the following expression for the {\Cp} fractional derivative, 
\begin{IEEEeqnarray}{l}\label{bruno1}
D_x^{\alpha}[f(g(x))] = 
\sum^{\infty}_{k= 0} 
\frac{ \sin [\pi  (\alpha-k)]}{\pi  (\alpha-k)} \left(\frac{ \Gamma (\alpha+1)}{\Gamma (k+1)} 
\right)
x^{k-\alpha}
\times \\\nn
\times
\sum _{m=0}^{k} \frac{1}{m!}\left(\sum _{j=0}^m (-1)^j \
{C^{j}_{m}}
g(x)^j
\frac{d^{k}}{dx^{k}}
g(x)^{m-j}\right) f^{(m)}(g(x))
.
\end{IEEEeqnarray}
As in the case of the fractional product rule, we exchange the summation order,
\begin{eqnarray}
\sum^{\infty}_{k= 0} \sum _{m=0}^{k} =
\sum _{m=0}^{\infty} \sum^{\infty}_{k=m}
,
\end{eqnarray}
which directly leads to, 
\begin{IEEEeqnarray}{l}\label{bruno2}
D_x^{\alpha}[f(g(x))] = 
\sum _{m=0}^{\infty}
\frac{f^{(m)}[g(x)]}{m!}
\sum^{\infty}_{k= m} 
\frac{ \sin [\pi  (\alpha-k)]}{\pi  (\alpha-k)} \left(\frac{ \Gamma (\alpha+1)}{\Gamma (k+1)} 
\right)
x^{k-\alpha} 
\times \\\nn
\times
\left(\sum _{j=0}^m (-1)^j \
{C^{j}_{m}}
g(x)^j
\frac{d^{k}}{dx^{k}}
g(x)^{m-j}\right) = \\\nn = 
\sum _{m=0}^{\infty}
\frac{f^{(m)}[g(x)]}{m!}
\sum _{j=0}^m (-1)^j \
{C^{j}_{m}}
g(x)^j
\sum^{\infty}_{k= m} 
\frac{ \sin [\pi  (\alpha-k)]}{\pi  (\alpha-k)} \left(\frac{ \Gamma (\alpha+1)}{\Gamma (k+1)} 
\right)
x^{k-\alpha} 
\frac{d^{k}}{dx^{k}}
g(x)^{m-j}.
\end{IEEEeqnarray}
By shifting the summation index in the inner-most sum in Eq.(\ref{bruno2}) as 
$k\to{k+m}$, we can rewrite it as
\begin{IEEEeqnarray}{l}\label{bruno3}
S_{1}(x,\alpha)\equiv{}\sum^{\infty}_{k= m} 
\frac{ \sin [\pi  (\alpha-k)]}{\pi  (\alpha-k)} \left(\frac{ \Gamma (\alpha+1)}{\Gamma (k+1)} 
\right)
x^{k-\alpha} 
\frac{d^{k}}{dx^{k}}
g(x)^{m-j}= \\\nn = 
x^{m-\alpha}
\left[
\sum^{\infty}_{k=0} 
\frac{ \sin [\pi  (\alpha-k-m)]}{\pi  (\alpha-k-m)} \left(\frac{ \Gamma (\alpha+1)}{\Gamma (k+m+1)} 
\right)
x^{k} 
\frac{d^{k}}{dx^{k}}
\right]
\frac{d^{m}}{dx^{m}}
g(x)^{m-j}.
\end{IEEEeqnarray}
In order to perform the summation in Eq.(\ref{bruno3}) we shall introduce the 
Cauchy--Euler differential operator $\widehat{J}_{x}\equiv{x\;{d}/{dx}}$ \cite{herrmann2014fractional}, along with the normal ordering operation \cite{peskin1995introduction},
\begin{eqnarray}
\cl\widehat{J}_{x}^{k} \cl =\;\;\cl \left(x\frac{d}{dx}\right)^{k} \cl = 
x^{k}\frac{d^{k}}{dx^{k}}
.
\end{eqnarray}
As the result, we obtain,
\begin{IEEEeqnarray}{l}\label{bruno4}
S_{2}(x,\alpha)\equiv{}
\sum^{\infty}_{k=0} 
\frac{ \sin [\pi  (\alpha-k-m)]}{\pi  (\alpha-k-m)} \left(\frac{ \Gamma (\alpha+1)}{\Gamma (k+m+1)} 
\right)
\cl\widehat{J}_{x}^{k} \cl  = \\\nn = 
\frac{\sin [\pi  (\alpha - m)]}{\pi  (\alpha - m) }
\frac{\Gamma (\alpha+1) }{\Gamma (m+1)}
\, \cl_{2}F_{2} \cl
\left[ 
\begin{array}{cc}
1, & m-\alpha\\ 
1+m, & 1+ m-\alpha
\end{array} ;\;\;\; - \widehat{J}_{x}  \right],
\end{IEEEeqnarray}
where $\cl_{2}F_{2} \cl$ is the normal ordered generalized hypergeometric function. The summation result obtained in Eq.(\ref{bruno4}) allows us to bring the sum in Eq.(\ref{bruno3}) into the following form,
\begin{equation}
S_{1}(x,\alpha)= 
x^{m-\alpha}
\frac{\sin [\pi  (\alpha - m)]}{\pi  (\alpha - m) }
\frac{\Gamma (\alpha+1) }{\Gamma (m+1)}
\, \cl_{2}F_{2} \cl
\left[ 
\begin{array}{cc}
1, & m-\alpha\\ 
1+m, & 1+ m-\alpha
\end{array} ;\;\;\; - \widehat{J}_{x}  \right] 
\frac{d^{m}}{dx^{m}}
g(x)^{m-j}.
\end{equation}
Finally, by introducing the weight function as,
\begin{equation}\label{weight_func}
W_{m}(\alpha,x,g(x)) = 
\sum _{j=0}^m 
\frac{(-1)^j}{m!} \
{C^{j}_{m}}
g(x)^j
\, \cl_{2}F_{2} \cl
\left[ 
\begin{array}{cc}
1, & m-\alpha\\ 
1+m, & 1 + m-\alpha
\end{array} ; \;\;- x\frac{d}{dx}  \right] 
\frac{d^{m}}{dx^{m}}
g(x)^{m-j}
,
\end{equation}
the fractional chain rule given by Eq.(\ref{bruno2}) acquires a particularly simple form,
\begin{eqnarray}\label{frac_chain_rule1}
D_x^{\alpha}[f(g(x))] = 
\sum _{m=0}^{\infty}
W_{m}(\alpha,x,g(x)) 
\frac{\sin [\pi  (\alpha  - m)]}{\pi  (\alpha -m) }
\frac{\Gamma (\alpha+1) }{\Gamma (m+1)}
 x^{m-\alpha}f^{(m)}[g(x)].
\end{eqnarray}
If we compare the obtained result given by Eq.~(\ref{frac_chain_rule1}) with the 
regular expression of the fractional derivative given by Eq.~(\ref{frac_integer_sum1}), we find that the fractional chain rule effectively reduces to the regular fractional derivative of a function $f(x)$ with the extra weight factor $W(x)$ given by Eq.~(\ref{weight_func}).

\section{ Applications  of the fractional chain and product rules}\label{section_appl}

The goal of this section is to apply the obtained fractional chain and product rules
to the exact evaluation of the  {\Cp} fractional derivative of hyperbolic tangent. First, we shall point out that the derived fractional product rule given by Eq.~(\ref{caputo_product}) leads to an {\it implicit} evaluation of the {\Cp} fractional derivative of a product of two functions. Indeed, the {\Cp} fractional derivative of order $\alpha$ of a product of two functions is expressed in terms of a semi-infinite sum of a product of an integer derivative of the $l^{\rm th}$ order of the first function and the {\Cp} fractional derivative of the $(\alpha-l)^{\rm th}$ order of the second function. 
In the particular case of elementary functions, the {\Cp} fractional derivative of the $(\alpha-l)^{\rm th}$ order results in the {\it repeated integral} of the generalized hypergeometric function with a power-law argument. In this paper, we implement the 
generalized Euler's integral transform developed in \cite{shchedrin2017exact} that allows us to represent the implicit form of the repeated integral in terms of a single hypergeometric function of a higher order. In this section, we will evaluate the {\Cp} fractional derivative of a composite function, e.g., hyperbolic tangent, that represents 
the dark soliton solution to the nonlinear {\Sch} equation \cite{boyd2003nonlinear}.
The obtained result is especially important for the evaluation of the dark soliton's kinetic energy in the course of generalization from the integer to the {\it fractional} non-linear {\Sch} equation.

We start with the direct application of the fractional product rule given by Eq.~(\ref{caputo_product}) to hyperbolic tangent,
\begin{IEEEeqnarray}{l}
{}^{C}D_x^{\alpha}\left[\frac{\sinh(\beta x)}{\cosh(\beta x)}\right] =
\sum^{\infty}_{l= 0} 
C^{l}_{\alpha}
\frac{d^{l}}{dx^{l}}\left(\frac{1}{\cosh(\beta x)}\right)
D_x^{(\alpha-l)}[\sinh(\beta x)]
.
\end{IEEEeqnarray}
First, we evaluate the $l^{\rm th}$ order integer derivative of hyperbolic secant by means of 
di Bruno's formula for an inverse function  \cite{faa1855sullo, di1857note, jordan1965calculus, wolfram},
\begin{eqnarray}
\frac{d^{n}}{dx^{n}}
\left(\frac{1}{f(x)}\right)=
(n+1) \sum _{k=0}^n  C_n^k \frac{(-1)^k }{k+1} 
\frac{1}{f(x)^{k+1}}
\frac{d^{n}}{dx^{n}}
f(x)^{k}.
\end{eqnarray}
Next, we evaluate the {\Cp} fractional derivative of hyperbolic sine, which can be done exactly by means of the generalized Euler's integral transform  \cite{shchedrin2017exact},
\begin{eqnarray}\label{caputo_sinh1}
{}^{C}D_x^{(\alpha)}[\sinh(\beta x)] = 
\frac{\beta  x^{1-\alpha } \, _1F_2\left(1;\frac{2-\alpha }{2},\frac{3-\alpha }{2}; 
\frac{\beta ^2 x^{2}}{4}\right)}{\Gamma (2-\alpha )}
.
\end{eqnarray}
Now we rewrite the fractional derivative of the $(\alpha-l)^{\rm th}$ order as, 
\begin{eqnarray}{\label{caputo_rep}}
D_x^{(\alpha-l)}[\sinh(\beta x)] = 
D_x^{(-l)} D_x^{(\alpha)}[\sinh(\beta x)]
=
\underbrace{\int^{x}_{0} {dx}\cdots \int^{x}_{0} {dx} }_{l \; \text {times}}
D_x^{(\alpha)}[\sinh(\beta x)]
.
\end{eqnarray}
In order to evaluate the repeated integral in Eq.~(\ref{caputo_rep}), we employ the {\Cy} repeated integration formula \cite{samko1993fractional, herrmann2014fractional},
\begin{eqnarray}
I_{n} \equiv \underbrace{\int^{x}_{0} {dx}\cdots \int^{x}_{0} {dx} }_{n \; \text {times}}f(x) = 
\frac{1}{\Gamma(n)} \int^{x}_{0} {dt} \; (x-t)^{n-1} f(t) 
.
\end{eqnarray}
Next, we apply the {\Cy} formula to a generalized hypergeometric function ${}_{A}F_{B}$, followed by re-scaling $t\to{xt}$, 
\begin{IEEEeqnarray}{l}\label{hypergeometry1}
J_{n} \equiv \underbrace{\int^{x}_{0} {dx}\cdots \int^{x}_{0} {dx} }_{n \; \text {times}}
\; x^{\kappa}
{}_{A}F_{B}\left[ 
\begin{array}{c}
a_{1},\ldots ,a_{A} \\ 
b_{1},\ldots ,b_{B}
\end{array} ; \zeta  x^{m} \right] 
 = 
\frac{1}{\Gamma(n)}
\int^{x}_{0}{dt} (x-t)^{n-1}
t^{\kappa}
{}_{A}F_{B}\left[ 
\begin{array}{c}
a_{1},\ldots ,a_{A} \\ 
b_{1},\ldots ,b_{B}
\end{array} ; \zeta  t^{m} \right] 
\nn\\ = 
\frac{x^{\kappa+n}}{\Gamma(n)}
\int^{1}_{0}{dt}\; t^{\kappa} (1-t)^{n-1}
{}_{A}F_{B}\left[ 
\begin{array}{c}
a_{1},\ldots ,a_{A} \\ 
b_{1},\ldots ,b_{B}
\end{array} ; \zeta x^{m}  t^{m} \right],
\end{IEEEeqnarray}
Where $\kappa$ and $\zeta$ are arbitrary real constants and $m$ is an arbitrary integer. We immediately recognize that the obtained integral in Eq.~(\ref{hypergeometry1}) is nothing but the generalized Euler's integral transform \cite{shchedrin2017exact}, which 
can be written as,
\begin{IEEEeqnarray}{l}\label{genresult1}
\int_{0}^{1} {dt}\;  t^{c-1} (1-t)_{{}}^{d-c-1}\ 
{}_{A}F_{B}\left[ 
\begin{array}{c}
a_{1},\dotsb ,a_{A} \\ 
b_{1},\dotsb ,b_{B}
\end{array} ; z  t^{m} \right]  
=\\\nn = 
\frac{\Gamma(d-c)\Gamma(c)}{\Gamma(d)}
{}_{A+m}F_{B+m}\left[ 
\begin{array}{cc}
a_{1},\dotsb ,a_{A}, & c_{0},\cdots c_{m-1} \\ 
b_{1},\dotsb ,b_{B},  & d_{0},\cdots d_{m-1}
\end{array} ; z \right],
\end{IEEEeqnarray}
where the constants $c_{j}$ and $d_{j}$ are given by $c_{j} ={(c + j)}/{m}$, and $d_{j} = {(d + j)}/{m}$ with index $j$ spanning $j\in[0,1,\cdots, m-1]$. Therefore, by assigning $z=\zeta x^m$, the repeated integral of a generalized hypergeometric function that appears in Eq.(\ref{hypergeometry1}) can be expressed in terms of a single generalized hypergeometric function of a higher order,
\begin{IEEEeqnarray}{l}\label{eit_applied1}
J_{n} =  x^{\kappa+n}
\frac{\Gamma(\kappa+1)}{\Gamma(\kappa+n+1)}
{}_{A+m}F_{B+m}\left[ 
\begin{array}{cccc}
a_{1},\ldots ,a_{A}, & c_{1} & \cdots & c_{m} \\ 
b_{1},\ldots ,b_{B}, & d_{1} &\cdots & d_{m} 
\end{array} ; \zeta x^m \right]
,
\end{IEEEeqnarray}
where $c_{j}={(\kappa+j)}/{m}$ and $d_{j}=(\kappa+n+j)/m$ for $j\in[1,2,\cdots, m]$.
Thus, a straightforward application of the general result given by Eq.~(\ref{eit_applied1})
with specific values $\kappa=1-\alpha$, $\zeta=\beta^2/4$, and $m=2$ results in an explicit expression for the Caputo fractional derivative of the $(\alpha-l)^{\text{th}}$ order of hyperbolic sine 
presented in Eq.(\ref{caputo_sinh1}),
\begin{IEEEeqnarray}{l}
{}^{C}D_x^{(\alpha-l)}[\sinh(\beta x)] = 
{}^{C}D_x^{-l}\left[
\frac{\beta  x^{1-\alpha } \, _1F_2\left(1;\frac{2-\alpha }{2},\frac{3-\alpha }{2}; 
\frac{\beta ^2 x^{2}}{4}\right)}{\Gamma (2-\alpha )} \right]=\\\nn = 
\frac{\beta}{\Gamma(2-\alpha)}
x^{1-\alpha+l}
\frac{\Gamma(2-\alpha)}{\Gamma(2+l-\alpha)}\;
{}_{3}F_{4}\left[ 
\begin{array}{cccc}
1, & \frac{2-\alpha }{2}, \frac{3-\alpha }{2},\\ 
\frac{2-\alpha }{2},\frac{3-\alpha }{2} & ,\frac{2+l-\alpha }{2}, \frac{3+l-\alpha }{2}
\end{array} ; \frac{\beta ^2 x^{2}}{4} \right] =\\\nn = 
\frac{\beta  x^{1-\alpha+l}}
{\Gamma(2+l-\alpha)}
\; {}_{1}F_{2}\left(1;\frac{2+l-\alpha }{2},\frac{3+l-\alpha }{2}; 
\frac{\beta ^2 x^{2}}{4}\right)
.
\end{IEEEeqnarray}

 \begin{figure}[t!]
\centering
\hspace{-2mm}
{\includegraphics[width=0.7\columnwidth]{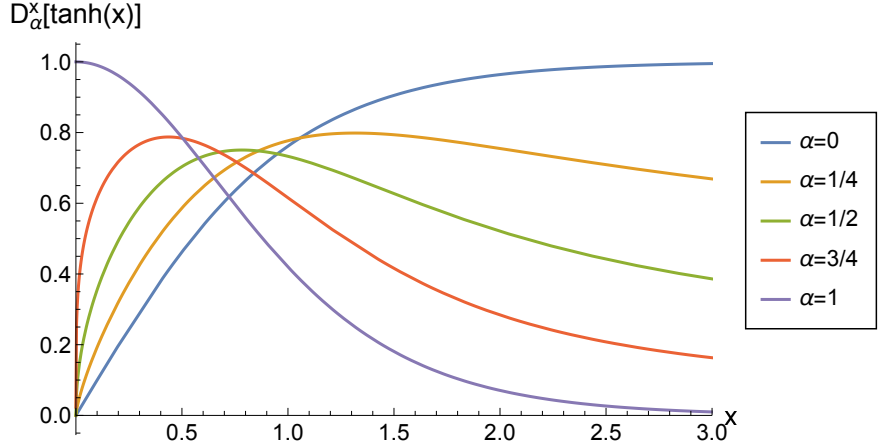}}
\caption{ The Caputo fractional derivative of $f(x)=\tanh (x)$ for fractional order $0\leq \alpha \leq 1$, with the specific choices of $\alpha$ displayed in the legend. This derivative was calculated using the expression for the fractional chain rule of the Caputo derivative, which contains an infinite sum of generalized hypergeometric functions. The series was truncated after 10 terms to produce the plots above.}
\label{fig2}
\end{figure}

Therefore, we obtain the sum of generalized hypergeometric functions for the fractional derivative of hyperbolic tangent,
\begin{IEEEeqnarray}{l}
{}^{C}D_x^{\alpha}\left[\frac{\sinh(\beta x)}{\cosh(\beta x)} \right] =
\sum^{\infty}_{l= 0} 
C^{l}_{\alpha}
\frac{d^{l}}{dx^{l}}\left(\frac{1}{\cosh(\beta x)}\right)
D_x^{(\alpha-l)}[\sinh(\beta x)]= \\\nn = 
\sum^{\infty}_{l= 0} 
C^{l}_{\alpha}
\frac{d^{l}}{dx^{l}}\left(\frac{1}{\cosh(\beta x)}\right)
\frac{\beta  x^{1-\alpha + l}}
{\Gamma(2+l-\alpha)}
\; {}_{1}F_{2}\left(1;\frac{2+l-\alpha }{2},\frac{3+l-\alpha }{2}; 
\frac{\beta ^2 x^{2}}{4}\right)
.
\end{IEEEeqnarray}
In a similar way, one can calculate the fractional derivative of products and ratios of elementary functions by expressing them in terms of generalized hypergeometric functions. The application of the {\Cy} formula to the repeated integral of a hypergeometric function, followed by the generalized EIT method, leads to an explicit result for the {\Cp} fractional derivative of an arbitrary function. However, in contrast to the generalized Euler integral transform that yields a single generalized hypergeometric function, both chain and product rules for the {\Cp}  fractional derivative produce an infinite series of generalized hypergeometric functions.

\section{Conclusions}\label{section_concl}

In conclusion, we explored the {\Cp} fractional derivative of composite functions, which, upon expansion into an infinite series of integer order derivatives, led to the derivation of the fractional product and chain rules. The practical implementation of these results was restricted since the final results were expressed in terms of a sum of the repeated integral of the generalized hypergeometric function. We applied the Euler's integral transform, which allowed us to transform the nested integral of generalized hypergeometric functions into a single hypergeometric function of a higher order.
Thus, we were able to obtain the exact result for the {\Cp} fractional derivative of a composite function, such as hyperbolic tangent. Moreover, fractional product and chain rules allowed us to extend the applicability of the generalized Euler's integral transform as a method for the exact evaluation of the {\Cp} fractional derivative of composite functions. However, unlike the Euler's integral transform that results in a single generalized hypergeometric function, both Caputo fractional chain and product rules produce an infinite series of hypergeometric functions.

\section*{Acknowledgments}
The authors would like to thank Daniel Jaschke and Marc Valdez for numerous and fruitful discussions. 
The authors acknowledge support from the US National Science Foundation under grant numbers PHY-1520915, OAC-1740130, DMR-1407962, and the US Air Force Office of Scientific Research grant number FA9550-14-1-0287. This work was performed in part at the Aspen Center for Physics, which is supported by the US National Science Foundation grant PHY-1607611.


\begin{thebibliography}{46}%
\makeatletter
\providecommand \@ifxundefined [1]{%
 \@ifx{#1\undefined}
}%
\providecommand \@ifnum [1]{%
 \ifnum #1\expandafter \@firstoftwo
 \else \expandafter \@secondoftwo
 \fi
}%
\providecommand \@ifx [1]{%
 \ifx #1\expandafter \@firstoftwo
 \else \expandafter \@secondoftwo
 \fi
}%
\providecommand \natexlab [1]{#1}%
\providecommand \enquote  [1]{``#1''}%
\providecommand \bibnamefont  [1]{#1}%
\providecommand \bibfnamefont [1]{#1}%
\providecommand \citenamefont [1]{#1}%
\providecommand \href@noop [0]{\@secondoftwo}%
\providecommand \href [0]{\begingroup \@sanitize@url \@href}%
\providecommand \@href[1]{\@@startlink{#1}\@@href}%
\providecommand \@@href[1]{\endgroup#1\@@endlink}%
\providecommand \@sanitize@url [0]{\catcode `\\12\catcode `\$12\catcode
  `\&12\catcode `\#12\catcode `\^12\catcode `\_12\catcode `\%12\relax}%
\providecommand \@@startlink[1]{}%
\providecommand \@@endlink[0]{}%
\providecommand \url  [0]{\begingroup\@sanitize@url \@url }%
\providecommand \@url [1]{\endgroup\@href {#1}{\urlprefix }}%
\providecommand \urlprefix  [0]{URL }%
\providecommand \Eprint [0]{\href }%
\providecommand \doibase [0]{http://dx.doi.org/}%
\providecommand \selectlanguage [0]{\@gobble}%
\providecommand \bibinfo  [0]{\@secondoftwo}%
\providecommand \bibfield  [0]{\@secondoftwo}%
\providecommand \translation [1]{[#1]}%
\providecommand \BibitemOpen [0]{}%
\providecommand \bibitemStop [0]{}%
\providecommand \bibitemNoStop [0]{.\EOS\space}%
\providecommand \EOS [0]{\spacefactor3000\relax}%
\providecommand \BibitemShut  [1]{\csname bibitem#1\endcsname}%
\let\auto@bib@innerbib\@empty
\bibitem [{\citenamefont {Murray}(1926)}]{murray1926physiological}%
  \BibitemOpen
  \bibfield  {author} {\bibinfo {author} {\bibfnamefont {Cecil~D}\ \bibnamefont
  {Murray}},\ }\bibfield  {title} {\enquote {\bibinfo {title} {The
  physiological principle of minimum work {I}. {T}he vascular system and the
  cost of blood volume},}\ }\href@noop {} {\bibfield  {journal} {\bibinfo
  {journal} {P. N. A. S.}\ }\textbf {\bibinfo {volume} {12}},\ \bibinfo {pages}
  {207} (\bibinfo {year} {1926})}\BibitemShut {NoStop}%
\bibitem [{\citenamefont {West}\ \emph {et~al.}(1997)\citenamefont {West},
  \citenamefont {Brown},\ and\ \citenamefont {Enquist}}]{west1997general}%
  \BibitemOpen
  \bibfield  {author} {\bibinfo {author} {\bibfnamefont {Geoffrey~B}\
  \bibnamefont {West}}, \bibinfo {author} {\bibfnamefont {James~H}\
  \bibnamefont {Brown}}, \ and\ \bibinfo {author} {\bibfnamefont {Brian~J}\
  \bibnamefont {Enquist}},\ }\bibfield  {title} {\enquote {\bibinfo {title} {A
  general model for the origin of allometric scaling laws in biology},}\
  }\href@noop {} {\bibfield  {journal} {\bibinfo  {journal} {Science}\ }\textbf
  {\bibinfo {volume} {276}},\ \bibinfo {pages} {122} (\bibinfo {year}
  {1997})}\BibitemShut {NoStop}%
\bibitem [{\citenamefont {Sherman}(1981)}]{sherman1981connecting}%
  \BibitemOpen
  \bibfield  {author} {\bibinfo {author} {\bibfnamefont {Thomas~F}\
  \bibnamefont {Sherman}},\ }\bibfield  {title} {\enquote {\bibinfo {title} {On
  connecting large vessels to small. {T}he meaning of {M}urray's law},}\
  }\href@noop {} {\bibfield  {journal} {\bibinfo  {journal} {J. Gen. Phys.}\
  }\textbf {\bibinfo {volume} {78}},\ \bibinfo {pages} {431} (\bibinfo {year}
  {1981})}\BibitemShut {NoStop}%
\bibitem [{\citenamefont {Lundstrom}\ \emph {et~al.}(2008)\citenamefont
  {Lundstrom}, \citenamefont {Higgs}, \citenamefont {Spain},\ and\
  \citenamefont {Fairhall}}]{lundstrom2008fractional}%
  \BibitemOpen
  \bibfield  {author} {\bibinfo {author} {\bibfnamefont {Brian~N}\ \bibnamefont
  {Lundstrom}}, \bibinfo {author} {\bibfnamefont {Matthew~H}\ \bibnamefont
  {Higgs}}, \bibinfo {author} {\bibfnamefont {William~J}\ \bibnamefont
  {Spain}}, \ and\ \bibinfo {author} {\bibfnamefont {Adrienne~L}\ \bibnamefont
  {Fairhall}},\ }\bibfield  {title} {\enquote {\bibinfo {title} {Fractional
  differentiation by neocortical pyramidal neurons},}\ }\href@noop {}
  {\bibfield  {journal} {\bibinfo  {journal} {Nat. Neuro.}\ }\textbf {\bibinfo
  {volume} {11}},\ \bibinfo {pages} {1335} (\bibinfo {year}
  {2008})}\BibitemShut {NoStop}%
\bibitem [{\citenamefont {Anastasio}(1994)}]{anastasio1994fractional}%
  \BibitemOpen
  \bibfield  {author} {\bibinfo {author} {\bibfnamefont {Thomas~J}\
  \bibnamefont {Anastasio}},\ }\bibfield  {title} {\enquote {\bibinfo {title}
  {The fractional-order dynamics of brainstem vestibulo-oculomotor neurons},}\
  }\href@noop {} {\bibfield  {journal} {\bibinfo  {journal} {Bio. Cyber.}\
  }\textbf {\bibinfo {volume} {72}},\ \bibinfo {pages} {69} (\bibinfo {year}
  {1994})}\BibitemShut {NoStop}%
\bibitem [{\citenamefont {Aizenberg}\ \emph {et~al.}(2005)\citenamefont
  {Aizenberg}, \citenamefont {Weaver}, \citenamefont {Thanawala}, \citenamefont
  {Sundar}, \citenamefont {Morse},\ and\ \citenamefont
  {Fratzl}}]{aizenberg2005skeleton}%
  \BibitemOpen
  \bibfield  {author} {\bibinfo {author} {\bibfnamefont {Joanna}\ \bibnamefont
  {Aizenberg}}, \bibinfo {author} {\bibfnamefont {James~C}\ \bibnamefont
  {Weaver}}, \bibinfo {author} {\bibfnamefont {Monica~S}\ \bibnamefont
  {Thanawala}}, \bibinfo {author} {\bibfnamefont {Vikram~C}\ \bibnamefont
  {Sundar}}, \bibinfo {author} {\bibfnamefont {Daniel~E}\ \bibnamefont
  {Morse}}, \ and\ \bibinfo {author} {\bibfnamefont {Peter}\ \bibnamefont
  {Fratzl}},\ }\bibfield  {title} {\enquote {\bibinfo {title} {Skeleton of
  euplectella sp.: structural hierarchy from the nanoscale to the
  macroscale},}\ }\href@noop {} {\bibfield  {journal} {\bibinfo  {journal}
  {Science}\ }\textbf {\bibinfo {volume} {309}},\ \bibinfo {pages} {275}
  (\bibinfo {year} {2005})}\BibitemShut {NoStop}%
\bibitem [{\citenamefont {Liu}\ and\ \citenamefont {Jiang}(2011)}]{liu2011bio}%
  \BibitemOpen
  \bibfield  {author} {\bibinfo {author} {\bibfnamefont {Kesong}\ \bibnamefont
  {Liu}}\ and\ \bibinfo {author} {\bibfnamefont {Lei}\ \bibnamefont {Jiang}},\
  }\bibfield  {title} {\enquote {\bibinfo {title} {Bio-inspired design of
  multiscale structures for function integration},}\ }\href@noop {} {\bibfield
  {journal} {\bibinfo  {journal} {Nano Today}\ }\textbf {\bibinfo {volume}
  {6}},\ \bibinfo {pages} {155} (\bibinfo {year} {2011})}\BibitemShut {NoStop}%
\bibitem [{\citenamefont {McCulloh}\ \emph {et~al.}(2003)\citenamefont
  {McCulloh}, \citenamefont {Sperry},\ and\ \citenamefont
  {Adler}}]{mcculloh2003water}%
  \BibitemOpen
  \bibfield  {author} {\bibinfo {author} {\bibfnamefont {Katherine~A}\
  \bibnamefont {McCulloh}}, \bibinfo {author} {\bibfnamefont {John~S}\
  \bibnamefont {Sperry}}, \ and\ \bibinfo {author} {\bibfnamefont
  {Frederick~R}\ \bibnamefont {Adler}},\ }\bibfield  {title} {\enquote
  {\bibinfo {title} {Water transport in plants obeys {M}urray's law},}\
  }\href@noop {} {\bibfield  {journal} {\bibinfo  {journal} {Nature}\ }\textbf
  {\bibinfo {volume} {421}},\ \bibinfo {pages} {939} (\bibinfo {year}
  {2003})}\BibitemShut {NoStop}%
\bibitem [{\citenamefont {Metzler}\ and\ \citenamefont
  {Klafter}(2000)}]{metzler2000random}%
  \BibitemOpen
  \bibfield  {author} {\bibinfo {author} {\bibfnamefont {Ralf}\ \bibnamefont
  {Metzler}}\ and\ \bibinfo {author} {\bibfnamefont {Joseph}\ \bibnamefont
  {Klafter}},\ }\bibfield  {title} {\enquote {\bibinfo {title} {The random
  walk's guide to anomalous diffusion: a fractional dynamics approach},}\
  }\href@noop {} {\bibfield  {journal} {\bibinfo  {journal} {Phys. Rep.}\
  }\textbf {\bibinfo {volume} {339}},\ \bibinfo {pages} {1} (\bibinfo {year}
  {2000})}\BibitemShut {NoStop}%
\bibitem [{\citenamefont {Saichev}\ and\ \citenamefont
  {Zaslavsky}(1997)}]{saichev1997fractional}%
  \BibitemOpen
  \bibfield  {author} {\bibinfo {author} {\bibfnamefont {Alexander~I}\
  \bibnamefont {Saichev}}\ and\ \bibinfo {author} {\bibfnamefont {George~M}\
  \bibnamefont {Zaslavsky}},\ }\bibfield  {title} {\enquote {\bibinfo {title}
  {Fractional kinetic equations: solutions and applications},}\ }\href@noop {}
  {\bibfield  {journal} {\bibinfo  {journal} {Chaos: An Interdisciplinary
  Journal of Nonlinear Science}\ }\textbf {\bibinfo {volume} {7}},\ \bibinfo
  {pages} {753} (\bibinfo {year} {1997})}\BibitemShut {NoStop}%
\bibitem [{\citenamefont {Mathai}\ \emph {et~al.}(2009)\citenamefont {Mathai},
  \citenamefont {Saxena},\ and\ \citenamefont {Haubold}}]{mathai2009h}%
  \BibitemOpen
  \bibfield  {author} {\bibinfo {author} {\bibfnamefont {Arakaparampil~M}\
  \bibnamefont {Mathai}}, \bibinfo {author} {\bibfnamefont {Ram~Kishore}\
  \bibnamefont {Saxena}}, \ and\ \bibinfo {author} {\bibfnamefont {Hans~J}\
  \bibnamefont {Haubold}},\ }\href@noop {} {\emph {\bibinfo {title} {The
  {H}-function: theory and applications}}}\ (\bibinfo  {publisher} {Springer},\
  \bibinfo {year} {2009})\BibitemShut {NoStop}%
\bibitem [{\citenamefont {Samko}\ \emph {et~al.}(1993)\citenamefont {Samko},
  \citenamefont {Kilbas},\ and\ \citenamefont
  {Marichev}}]{samko1993fractional}%
  \BibitemOpen
  \bibfield  {author} {\bibinfo {author} {\bibfnamefont {Stefan~G}\
  \bibnamefont {Samko}}, \bibinfo {author} {\bibfnamefont {Anatoly~A}\
  \bibnamefont {Kilbas}}, \ and\ \bibinfo {author} {\bibfnamefont {Oleg~I}\
  \bibnamefont {Marichev}},\ }\href@noop {} {\emph {\bibinfo {title}
  {Fractional integrals and derivatives: Theory and Applications}}}\ (\bibinfo
  {publisher} {Gordon and Breach},\ \bibinfo {year} {1993})\BibitemShut
  {NoStop}%
\bibitem [{\citenamefont {Kilbas}\ \emph {et~al.}(2006)\citenamefont {Kilbas},
  \citenamefont {Srivastava},\ and\ \citenamefont
  {Trujillo}}]{kilbas2006theory}%
  \BibitemOpen
  \bibfield  {author} {\bibinfo {author} {\bibfnamefont {Anatoly~A}\
  \bibnamefont {Kilbas}}, \bibinfo {author} {\bibfnamefont {Hari~M}\
  \bibnamefont {Srivastava}}, \ and\ \bibinfo {author} {\bibfnamefont {Juan~J}\
  \bibnamefont {Trujillo}},\ }\href@noop {} {\emph {\bibinfo {title} {Theory
  and Applications of Fractional Differential Equations}}}\ (\bibinfo
  {publisher} {Elsevier},\ \bibinfo {year} {2006})\BibitemShut {NoStop}%
\bibitem [{\citenamefont {McBride}(1979)}]{mcbride1979fractional}%
  \BibitemOpen
  \bibfield  {author} {\bibinfo {author} {\bibfnamefont {Adam~C}\ \bibnamefont
  {McBride}},\ }\href@noop {} {\emph {\bibinfo {title} {Fractional calculus and
  integral transforms of generalized functions}}},\ Vol.~\bibinfo {volume}
  {31}\ (\bibinfo  {publisher} {Pitman London},\ \bibinfo {year}
  {1979})\BibitemShut {NoStop}%
\bibitem [{\citenamefont {West}(2014)}]{west2014colloquium}%
  \BibitemOpen
  \bibfield  {author} {\bibinfo {author} {\bibfnamefont {Bruce~J}\ \bibnamefont
  {West}},\ }\bibfield  {title} {\enquote {\bibinfo {title} {Colloquium:
  Fractional calculus view of complexity: A tutorial},}\ }\href@noop {}
  {\bibfield  {journal} {\bibinfo  {journal} {Rev. Mod. Phys.}\ }\textbf
  {\bibinfo {volume} {86}},\ \bibinfo {pages} {1169} (\bibinfo {year}
  {2014})}\BibitemShut {NoStop}%
\bibitem [{\citenamefont {Herrmann}(2014)}]{herrmann2014fractional}%
  \BibitemOpen
  \bibfield  {author} {\bibinfo {author} {\bibfnamefont {Richard}\ \bibnamefont
  {Herrmann}},\ }\href@noop {} {\emph {\bibinfo {title} {Fractional calculus:
  an introduction for physicists}}}\ (\bibinfo  {publisher} {World
  Scientific},\ \bibinfo {year} {2014})\BibitemShut {NoStop}%
\bibitem [{\citenamefont {Kiryakova}(1993)}]{kiryakova1993generalized}%
  \BibitemOpen
  \bibfield  {author} {\bibinfo {author} {\bibfnamefont {Virginia~S}\
  \bibnamefont {Kiryakova}},\ }\href@noop {} {\emph {\bibinfo {title}
  {Generalized fractional calculus and applications}}}\ (\bibinfo  {publisher}
  {CRC press},\ \bibinfo {year} {1993})\BibitemShut {NoStop}%
\bibitem [{\citenamefont {Zheng}\ \emph {et~al.}(2017)\citenamefont {Zheng},
  \citenamefont {Shen}, \citenamefont {Wang}, \citenamefont {Li}, \citenamefont
  {Dunphy}, \citenamefont {Hasan}, \citenamefont {Brinker},\ and\ \citenamefont
  {Su}}]{zheng2017bio}%
  \BibitemOpen
  \bibfield  {author} {\bibinfo {author} {\bibfnamefont {Xianfeng}\
  \bibnamefont {Zheng}}, \bibinfo {author} {\bibfnamefont {Guofang}\
  \bibnamefont {Shen}}, \bibinfo {author} {\bibfnamefont {Chao}\ \bibnamefont
  {Wang}}, \bibinfo {author} {\bibfnamefont {Yu}~\bibnamefont {Li}}, \bibinfo
  {author} {\bibfnamefont {Darren}\ \bibnamefont {Dunphy}}, \bibinfo {author}
  {\bibfnamefont {Tawfique}\ \bibnamefont {Hasan}}, \bibinfo {author}
  {\bibfnamefont {C~Jeffrey}\ \bibnamefont {Brinker}}, \ and\ \bibinfo {author}
  {\bibfnamefont {Bao-Lian}\ \bibnamefont {Su}},\ }\bibfield  {title} {\enquote
  {\bibinfo {title} {Bio-inspired {M}urray materials for mass transfer and
  activity},}\ }\href@noop {} {\bibfield  {journal} {\bibinfo  {journal}
  {Nature communications}\ }\textbf {\bibinfo {volume} {8}},\ \bibinfo {pages}
  {14921} (\bibinfo {year} {2017})}\BibitemShut {NoStop}%
\bibitem [{\citenamefont {Kilbas}\ \emph {et~al.}(2004)\citenamefont {Kilbas},
  \citenamefont {Saigo},\ and\ \citenamefont {Saxena}}]{kilbas2004generalized}%
  \BibitemOpen
  \bibfield  {author} {\bibinfo {author} {\bibfnamefont {Anatoly~A}\
  \bibnamefont {Kilbas}}, \bibinfo {author} {\bibfnamefont {Megumi}\
  \bibnamefont {Saigo}}, \ and\ \bibinfo {author} {\bibfnamefont
  {RK}~\bibnamefont {Saxena}},\ }\bibfield  {title} {\enquote {\bibinfo {title}
  {Generalized {M}ittag-{L}effler function and generalized fractional calculus
  operators},}\ }\href@noop {} {\bibfield  {journal} {\bibinfo  {journal} {Int.
  Trans. Spec. Fun.}\ }\textbf {\bibinfo {volume} {15}},\ \bibinfo {pages} {31}
  (\bibinfo {year} {2004})}\BibitemShut {NoStop}%
\bibitem [{\citenamefont {Srivastava}\ and\ \citenamefont
  {Tomovski}(2009)}]{srivastava2009fractional}%
  \BibitemOpen
  \bibfield  {author} {\bibinfo {author} {\bibfnamefont {HM}~\bibnamefont
  {Srivastava}}\ and\ \bibinfo {author} {\bibfnamefont {{\v{Z}}ivorad}\
  \bibnamefont {Tomovski}},\ }\bibfield  {title} {\enquote {\bibinfo {title}
  {Fractional calculus with an integral operator containing a generalized
  {M}ittag--{L}effler function in the kernel},}\ }\href@noop {} {\bibfield
  {journal} {\bibinfo  {journal} {Appl.Math. Comp.}\ }\textbf {\bibinfo
  {volume} {211}},\ \bibinfo {pages} {198} (\bibinfo {year}
  {2009})}\BibitemShut {NoStop}%
\bibitem [{\citenamefont {Shukla}\ and\ \citenamefont
  {Prajapati}(2007)}]{shukla2007generalization}%
  \BibitemOpen
  \bibfield  {author} {\bibinfo {author} {\bibfnamefont {AK}~\bibnamefont
  {Shukla}}\ and\ \bibinfo {author} {\bibfnamefont {JC}~\bibnamefont
  {Prajapati}},\ }\bibfield  {title} {\enquote {\bibinfo {title} {On a
  generalization of {M}ittag-{L}effler function and its properties},}\
  }\href@noop {} {\bibfield  {journal} {\bibinfo  {journal} {J. Math. Appl.}\
  }\textbf {\bibinfo {volume} {336}},\ \bibinfo {pages} {797} (\bibinfo {year}
  {2007})}\BibitemShut {NoStop}%
\bibitem [{\citenamefont {Mainardi}\ and\ \citenamefont
  {Gorenflo}(2000)}]{mainardi2000mittag}%
  \BibitemOpen
  \bibfield  {author} {\bibinfo {author} {\bibfnamefont {Francesco}\
  \bibnamefont {Mainardi}}\ and\ \bibinfo {author} {\bibfnamefont {Rudolf}\
  \bibnamefont {Gorenflo}},\ }\bibfield  {title} {\enquote {\bibinfo {title}
  {On {M}ittag-{L}effler-type functions in fractional evolution processes},}\
  }\href@noop {} {\bibfield  {journal} {\bibinfo  {journal} {J. Comp. Appl.
  Math.}\ }\textbf {\bibinfo {volume} {118}},\ \bibinfo {pages} {283} (\bibinfo
  {year} {2000})}\BibitemShut {NoStop}%
\bibitem [{\citenamefont {Haubold}\ \emph {et~al.}(2011)\citenamefont
  {Haubold}, \citenamefont {Mathai},\ and\ \citenamefont
  {Saxena}}]{haubold2011mittag}%
  \BibitemOpen
  \bibfield  {author} {\bibinfo {author} {\bibfnamefont {Hans~J}\ \bibnamefont
  {Haubold}}, \bibinfo {author} {\bibfnamefont {Arak~M}\ \bibnamefont
  {Mathai}}, \ and\ \bibinfo {author} {\bibfnamefont {Ram~K}\ \bibnamefont
  {Saxena}},\ }\bibfield  {title} {\enquote {\bibinfo {title}
  {{M}ittag-{L}effler functions and their applications},}\ }\href@noop {}
  {\bibfield  {journal} {\bibinfo  {journal} {J. Appl. Math.}\ }\textbf
  {\bibinfo {volume} {2011}},\ \bibinfo {pages} {1} (\bibinfo {year}
  {2011})}\BibitemShut {NoStop}%
\bibitem [{\citenamefont {Saadatmandi}\ and\ \citenamefont
  {Dehghan}(2010)}]{saadatmandi2010new}%
  \BibitemOpen
  \bibfield  {author} {\bibinfo {author} {\bibfnamefont {Abbas}\ \bibnamefont
  {Saadatmandi}}\ and\ \bibinfo {author} {\bibfnamefont {Mehdi}\ \bibnamefont
  {Dehghan}},\ }\bibfield  {title} {\enquote {\bibinfo {title} {A new
  operational matrix for solving fractional-order differential equations},}\
  }\href@noop {} {\bibfield  {journal} {\bibinfo  {journal} {Comp. \& Math
  Apps.}\ }\textbf {\bibinfo {volume} {59}},\ \bibinfo {pages} {1326} (\bibinfo
  {year} {2010})}\BibitemShut {NoStop}%
\bibitem [{\citenamefont {Shawagfeh}(2002)}]{shawagfeh2002analytical}%
  \BibitemOpen
  \bibfield  {author} {\bibinfo {author} {\bibfnamefont {Nabil~T}\ \bibnamefont
  {Shawagfeh}},\ }\bibfield  {title} {\enquote {\bibinfo {title} {Analytical
  approximate solutions for nonlinear fractional differential equations},}\
  }\href@noop {} {\bibfield  {journal} {\bibinfo  {journal} {Appl. Math.
  Comp.}\ }\textbf {\bibinfo {volume} {131}},\ \bibinfo {pages} {517} (\bibinfo
  {year} {2002})}\BibitemShut {NoStop}%
\bibitem [{\citenamefont {Momani}\ and\ \citenamefont
  {Odibat}(2006)}]{momani2006analytical}%
  \BibitemOpen
  \bibfield  {author} {\bibinfo {author} {\bibfnamefont {Shaher}\ \bibnamefont
  {Momani}}\ and\ \bibinfo {author} {\bibfnamefont {Zaid}\ \bibnamefont
  {Odibat}},\ }\bibfield  {title} {\enquote {\bibinfo {title} {Analytical
  solution of a time-fractional {N}avier--{S}tokes equation by {A}domian
  decomposition method},}\ }\href@noop {} {\bibfield  {journal} {\bibinfo
  {journal} {Appl. Math. Comp.}\ }\textbf {\bibinfo {volume} {177}},\ \bibinfo
  {pages} {488} (\bibinfo {year} {2006})}\BibitemShut {NoStop}%
\bibitem [{\citenamefont {Odibat}\ and\ \citenamefont
  {Momani}(2006)}]{odibat2006application}%
  \BibitemOpen
  \bibfield  {author} {\bibinfo {author} {\bibfnamefont {ZM}~\bibnamefont
  {Odibat}}\ and\ \bibinfo {author} {\bibfnamefont {Shaher}\ \bibnamefont
  {Momani}},\ }\bibfield  {title} {\enquote {\bibinfo {title} {Application of
  variational iteration method to nonlinear differential equations of
  fractional order},}\ }\href@noop {} {\bibfield  {journal} {\bibinfo
  {journal} {Int. J. Nonlin. Sc. and Num. Sim.}\ }\textbf {\bibinfo {volume}
  {7}},\ \bibinfo {pages} {27} (\bibinfo {year} {2006})}\BibitemShut {NoStop}%
\bibitem [{\citenamefont {Rossikhin}\ and\ \citenamefont
  {Shitikova}(2010)}]{rossikhin2010application}%
  \BibitemOpen
  \bibfield  {author} {\bibinfo {author} {\bibfnamefont {Yuriy~A}\ \bibnamefont
  {Rossikhin}}\ and\ \bibinfo {author} {\bibfnamefont {Marina~V}\ \bibnamefont
  {Shitikova}},\ }\bibfield  {title} {\enquote {\bibinfo {title} {Application
  of fractional calculus for dynamic problems of solid mechanics: novel trends
  and recent results},}\ }\href@noop {} {\bibfield  {journal} {\bibinfo
  {journal} {Appl. Mech. Rev.}\ }\textbf {\bibinfo {volume} {63}},\ \bibinfo
  {pages} {010801} (\bibinfo {year} {2010})}\BibitemShut {NoStop}%
\bibitem [{\citenamefont {Daftardar-Gejji}\ and\ \citenamefont
  {Jafari}(2005)}]{daftardar2005adomian}%
  \BibitemOpen
  \bibfield  {author} {\bibinfo {author} {\bibfnamefont {Varsha}\ \bibnamefont
  {Daftardar-Gejji}}\ and\ \bibinfo {author} {\bibfnamefont {Hossein}\
  \bibnamefont {Jafari}},\ }\bibfield  {title} {\enquote {\bibinfo {title}
  {Adomian decomposition: a tool for solving a system of fractional
  differential equations},}\ }\href@noop {} {\bibfield  {journal} {\bibinfo
  {journal} {J. Math. An. Appl.}\ }\textbf {\bibinfo {volume} {301}},\ \bibinfo
  {pages} {508} (\bibinfo {year} {2005})}\BibitemShut {NoStop}%
\bibitem [{\citenamefont {Dehghan}\ \emph {et~al.}(2010)\citenamefont
  {Dehghan}, \citenamefont {Manafian},\ and\ \citenamefont
  {Saadatmandi}}]{dehghan2010solving}%
  \BibitemOpen
  \bibfield  {author} {\bibinfo {author} {\bibfnamefont {Mehdi}\ \bibnamefont
  {Dehghan}}, \bibinfo {author} {\bibfnamefont {Jalil}\ \bibnamefont
  {Manafian}}, \ and\ \bibinfo {author} {\bibfnamefont {Abbas}\ \bibnamefont
  {Saadatmandi}},\ }\bibfield  {title} {\enquote {\bibinfo {title} {Solving
  nonlinear fractional partial differential equations using the homotopy
  analysis method},}\ }\href@noop {} {\bibfield  {journal} {\bibinfo  {journal}
  {Num. Meth. P. D. E.}\ }\textbf {\bibinfo {volume} {26}},\ \bibinfo {pages}
  {448} (\bibinfo {year} {2010})}\BibitemShut {NoStop}%
\bibitem [{\citenamefont {Momani}(2005)}]{momani2005explicit}%
  \BibitemOpen
  \bibfield  {author} {\bibinfo {author} {\bibfnamefont {Shaher}\ \bibnamefont
  {Momani}},\ }\bibfield  {title} {\enquote {\bibinfo {title} {An explicit and
  numerical solutions of the fractional {KdV} equation},}\ }\href@noop {}
  {\bibfield  {journal} {\bibinfo  {journal} {Math. Comp. Sim.}\ }\textbf
  {\bibinfo {volume} {70}},\ \bibinfo {pages} {110} (\bibinfo {year}
  {2005})}\BibitemShut {NoStop}%
\bibitem [{\citenamefont {Li}\ \emph {et~al.}(2010)\citenamefont {Li},
  \citenamefont {Chen},\ and\ \citenamefont {Podlubny}}]{li2010stability}%
  \BibitemOpen
  \bibfield  {author} {\bibinfo {author} {\bibfnamefont {Yan}\ \bibnamefont
  {Li}}, \bibinfo {author} {\bibfnamefont {YangQuan}\ \bibnamefont {Chen}}, \
  and\ \bibinfo {author} {\bibfnamefont {Igor}\ \bibnamefont {Podlubny}},\
  }\bibfield  {title} {\enquote {\bibinfo {title} {Stability of
  fractional-order nonlinear dynamic systems: {L}yapunov direct method and
  generalized {M}ittag--{L}effler stability},}\ }\href@noop {} {\bibfield
  {journal} {\bibinfo  {journal} {Comp.\& Math. Appl.}\ }\textbf {\bibinfo
  {volume} {59}},\ \bibinfo {pages} {1810} (\bibinfo {year}
  {2010})}\BibitemShut {NoStop}%
\bibitem [{\citenamefont {Diethelm}\ \emph {et~al.}(2005)\citenamefont
  {Diethelm}, \citenamefont {Ford}, \citenamefont {Freed},\ and\ \citenamefont
  {Luchko}}]{diethelm2005algorithms}%
  \BibitemOpen
  \bibfield  {author} {\bibinfo {author} {\bibfnamefont {Kai}\ \bibnamefont
  {Diethelm}}, \bibinfo {author} {\bibfnamefont {Neville~J}\ \bibnamefont
  {Ford}}, \bibinfo {author} {\bibfnamefont {Alan~D}\ \bibnamefont {Freed}}, \
  and\ \bibinfo {author} {\bibfnamefont {Yu}~\bibnamefont {Luchko}},\
  }\bibfield  {title} {\enquote {\bibinfo {title} {Algorithms for the
  fractional calculus: a selection of numerical methods},}\ }\href@noop {}
  {\bibfield  {journal} {\bibinfo  {journal} {Comp. Meth. Appl. Mech. Eng.}\
  }\textbf {\bibinfo {volume} {194}},\ \bibinfo {pages} {743} (\bibinfo {year}
  {2005})}\BibitemShut {NoStop}%
\bibitem [{\citenamefont {Agrawal}(2002)}]{agrawal2002solution}%
  \BibitemOpen
  \bibfield  {author} {\bibinfo {author} {\bibfnamefont {Om~P}\ \bibnamefont
  {Agrawal}},\ }\bibfield  {title} {\enquote {\bibinfo {title} {Solution for a
  fractional diffusion-wave equation defined in a bounded domain},}\
  }\href@noop {} {\bibfield  {journal} {\bibinfo  {journal} {Nonlin. Dyn.}\
  }\textbf {\bibinfo {volume} {29}},\ \bibinfo {pages} {145} (\bibinfo {year}
  {2002})}\BibitemShut {NoStop}%
\bibitem [{\citenamefont {Kiryakova}(2010)}]{kiryakova2010special}%
  \BibitemOpen
  \bibfield  {author} {\bibinfo {author} {\bibfnamefont {Virginia}\
  \bibnamefont {Kiryakova}},\ }\bibfield  {title} {\enquote {\bibinfo {title}
  {The special functions of fractional calculus as generalized fractional
  calculus operators of some basic functions},}\ }\href@noop {} {\bibfield
  {journal} {\bibinfo  {journal} {Comp. \& Math. Appl.}\ }\textbf {\bibinfo
  {volume} {59}},\ \bibinfo {pages} {1128} (\bibinfo {year}
  {2010})}\BibitemShut {NoStop}%
\bibitem [{\citenamefont {Kiryakova}(1997)}]{kiryakova1997all}%
  \BibitemOpen
  \bibfield  {author} {\bibinfo {author} {\bibfnamefont {Virginia}\
  \bibnamefont {Kiryakova}},\ }\bibfield  {title} {\enquote {\bibinfo {title}
  {All the special functions are fractional differintegrals of elementary
  functions},}\ }\href@noop {} {\bibfield  {journal} {\bibinfo  {journal} {J.
  Phys. A}\ }\textbf {\bibinfo {volume} {30}},\ \bibinfo {pages} {5085}
  (\bibinfo {year} {1997})}\BibitemShut {NoStop}%
\bibitem [{\citenamefont {Gladkina}\ \emph {et~al.}(2017)\citenamefont
  {Gladkina}, \citenamefont {Shchedrin}, \citenamefont {Khawaja},\ and\
  \citenamefont {Carr}}]{gladkina2017expansion}%
  \BibitemOpen
  \bibfield  {author} {\bibinfo {author} {\bibfnamefont {Anastasia}\
  \bibnamefont {Gladkina}}, \bibinfo {author} {\bibfnamefont {Gavriil}\
  \bibnamefont {Shchedrin}}, \bibinfo {author} {\bibfnamefont {U~Al}\
  \bibnamefont {Khawaja}}, \ and\ \bibinfo {author} {\bibfnamefont {Lincoln~D}\
  \bibnamefont {Carr}},\ }\bibfield  {title} {\enquote {\bibinfo {title}
  {Expansion of fractional derivatives in terms of an integer derivative
  series: physical and numerical applications},}\ }\href@noop {} {\bibfield
  {journal} {\bibinfo  {journal} {arXiv preprint arXiv:1710.06297}\ } (\bibinfo
  {year} {2017})}\BibitemShut {NoStop}%
\bibitem [{\citenamefont {Shchedrin}\ \emph {et~al.}(2017)\citenamefont
  {Shchedrin}, \citenamefont {Smith}, \citenamefont {Gladkina},\ and\
  \citenamefont {Carr}}]{shchedrin2017exact}%
  \BibitemOpen
  \bibfield  {author} {\bibinfo {author} {\bibfnamefont {Gavriil}\ \bibnamefont
  {Shchedrin}}, \bibinfo {author} {\bibfnamefont {Nathanael~C}\ \bibnamefont
  {Smith}}, \bibinfo {author} {\bibfnamefont {Anastasia}\ \bibnamefont
  {Gladkina}}, \ and\ \bibinfo {author} {\bibfnamefont {Lincoln~D}\
  \bibnamefont {Carr}},\ }\bibfield  {title} {\enquote {\bibinfo {title} {Exact
  results for a fractional derivative of elementary functions},}\ }\href@noop
  {} {\bibfield  {journal} {\bibinfo  {journal} {arXiv preprint
  arXiv:1711.07126}\ } (\bibinfo {year} {2017})}\BibitemShut {NoStop}%
\bibitem [{\citenamefont {Laskin}(2002)}]{laskin2002fractional}%
  \BibitemOpen
  \bibfield  {author} {\bibinfo {author} {\bibfnamefont {Nick}\ \bibnamefont
  {Laskin}},\ }\bibfield  {title} {\enquote {\bibinfo {title} {Fractional
  schr{\"o}dinger equation},}\ }\href@noop {} {\bibfield  {journal} {\bibinfo
  {journal} {Phys. Rev. E}\ }\textbf {\bibinfo {volume} {66}},\ \bibinfo
  {pages} {056108} (\bibinfo {year} {2002})}\BibitemShut {NoStop}%
\bibitem [{\citenamefont {Landau}\ and\ \citenamefont
  {Lifshits}(1965)}]{landau1965quantum}%
  \BibitemOpen
  \bibfield  {author} {\bibinfo {author} {\bibfnamefont {Lev~Davidovich}\
  \bibnamefont {Landau}}\ and\ \bibinfo {author} {\bibfnamefont
  {Evgenii~Mikhailovich}\ \bibnamefont {Lifshits}},\ }\href@noop {} {\emph
  {\bibinfo {title} {Quantum Mechanics: Non-relativistic Theory}}}\ (\bibinfo
  {publisher} {Pergamon Press},\ \bibinfo {year} {1965})\BibitemShut {NoStop}%
\bibitem [{\citenamefont {Faa~di Bruno}(1855)}]{faa1855sullo}%
  \BibitemOpen
  \bibfield  {author} {\bibinfo {author} {\bibfnamefont {Francesco}\
  \bibnamefont {Faa~di Bruno}},\ }\bibfield  {title} {\enquote {\bibinfo
  {title} {Sullo sviluppo delle funzioni},}\ }\href@noop {} {\bibfield
  {journal} {\bibinfo  {journal} {Annali di scienze matematiche e fisiche}\
  }\textbf {\bibinfo {volume} {6}},\ \bibinfo {pages} {479} (\bibinfo {year}
  {1855})}\BibitemShut {NoStop}%
\bibitem [{\citenamefont {Faa~di Bruno}(1857)}]{di1857note}%
  \BibitemOpen
  \bibfield  {author} {\bibinfo {author} {\bibfnamefont {Francesco}\
  \bibnamefont {Faa~di Bruno}},\ }\bibfield  {title} {\enquote {\bibinfo
  {title} {Note sur une nouvelle formule de calcul diff{\'e}rentiel},}\
  }\href@noop {} {\bibfield  {journal} {\bibinfo  {journal} {Quart. Jour. of
  Pure and App. Math}\ }\textbf {\bibinfo {volume} {1}},\ \bibinfo {pages}
  {359} (\bibinfo {year} {1857})}\BibitemShut {NoStop}%
\bibitem [{\citenamefont {Jord{\'a}n}(1965)}]{jordan1965calculus}%
  \BibitemOpen
  \bibfield  {author} {\bibinfo {author} {\bibfnamefont {K{\'a}roly}\
  \bibnamefont {Jord{\'a}n}},\ }\href@noop {} {\emph {\bibinfo {title}
  {Calculus of finite differences}}},\ Vol.~\bibinfo {volume} {33}\ (\bibinfo
  {publisher} {American Mathematical Soc.},\ \bibinfo {year}
  {1965})\BibitemShut {NoStop}%
\bibitem [{\citenamefont {Inc.}()}]{wolfram}%
  \BibitemOpen
  \bibfield  {author} {\bibinfo {author} {\bibfnamefont {Wolfram~Research{,}}\
  \bibnamefont {Inc.}},\ }\href
  {http://functions.wolfram.com/GeneralIdentities/9/} {\enquote {\bibinfo
  {title} {General mathematical identities for analytic functions:
  Differentiation},}\ }\bibinfo {note} {Champaign, IL, 2018}\BibitemShut
  {NoStop}%
\bibitem [{\citenamefont {Peskin}\ and\ \citenamefont
  {Schroeder}(1995)}]{peskin1995introduction}%
  \BibitemOpen
  \bibfield  {author} {\bibinfo {author} {\bibfnamefont {Michael~E}\
  \bibnamefont {Peskin}}\ and\ \bibinfo {author} {\bibfnamefont {Daniel~V}\
  \bibnamefont {Schroeder}},\ }\href@noop {} {\emph {\bibinfo {title} {An
  Introduction to Quantum Field Theory}}}\ (\bibinfo  {publisher} {Westview
  Press},\ \bibinfo {year} {1995})\BibitemShut {NoStop}%
\bibitem [{\citenamefont {Boyd}(2003)}]{boyd2003nonlinear}%
  \BibitemOpen
  \bibfield  {author} {\bibinfo {author} {\bibfnamefont {Robert~W}\
  \bibnamefont {Boyd}},\ }\href@noop {} {\emph {\bibinfo {title} {Nonlinear
  optics}}}\ (\bibinfo  {publisher} {Academic Press},\ \bibinfo {year}
  {2003})\BibitemShut {NoStop}%
\end{thebibliography}

%

\end{document}